\documentclass[preprint,12pt]{elsarticle}

\usepackage{amsmath,amssymb,amsthm}
\usepackage{enumerate}

\newtheorem{theorem}{Theorem}[section]
\newtheorem{lemma}{Lemma}[section]
\newtheorem{definition}{Definition}[section]
\newtheorem*{remark}{Remark}

\journal{Discrete Mathematics}

\begin{document}

\begin{frontmatter}



\title{A complete dichotomy theorem on the sparse $t$-Uniform Hypergraphicality Problem}

\author[label1,label2]{Istv\'an Mikl\'os}
\author[label1,label3]{Mikl\'os Ruszink\'o}
\author[label1]{Bogd\'an Zavalnij}

\affiliation[label1]{
  organization={HUN-REN Alfr\'ed R\'enyi Institute of Mathematics},
  addressline={Re\'altanoda u. 13--15},
  city={Budapest},
  postcode={1053},
  country={Hungary}
}

\affiliation[label2]{
  organization={HUN-REN SZTAKI},
  addressline={L\'agym\'anyosi u. 11},
  city={Budapest},
  postcode={1111},
  country={Hungary}
}

\affiliation[label3]{
  organization={Sorbonne University},
  addressline={Al Reem Island},
  city={Abu Dhabi},
  postcode={PO Box 38044},
  country={United Arab Emirates}
}



\begin{abstract}
We prove a complete dichotomy theorem for the parameterized sparse $t$-uniform hypergraphic degree sequence problem, $\mathrm{sparse}\text{-}t\text{-}\mathrm{uni}\text{-}\mathrm{HDS}_{\alpha',\alpha}$. For any fixed $t \ge 3$, given parameters $0 \le \alpha' \le \alpha < t-1$, the input consists of degree sequences $D$ of length $n$ with degrees between $n^{\alpha'}$ and $6n^{\alpha}$. We show that the problem is NP-complete whenever $\alpha' \le \frac{t(\alpha - 1) + 1}{t - 1}$, and solvable in linear time when $\alpha' > \frac{t(\alpha - 1) + 1}{t - 1}$. 
This establishes a sharp boundary between polynomial-time solvable and NP-complete instances,
thereby characterizing the computational complexity across all degree exponent regimes. The result extends the earlier NP-completeness of dense hypergraphicality to a unified framework covering both sparse and dense regimes, revealing that even extremely sparse instances (with maximum degree $o(n)$ but $\Omega(n^{\frac{t-1}{t}})$) remain NP-complete. On the other hand, the $t$-uniform hypergraphicality solvable in linear time when the maximum degree is $o(n^{\frac{t-1}{t}})$. This dichotomy  provides a comprehensive classification of the complexity landscape for hypergraphic degree sequences.
\end{abstract}




\end{frontmatter}



\section{Introduction}
The classical \emph{graphic degree sequence} problem asks whether a given sequence of integers can be realized as the vertex degrees of a simple graph.  Erd\H{o}s and Gallai~\cite{EG1960} gave necessary and sufficient conditions for a sequence to be graphic, and Havel~\cite{Havel1955} and Hakimi~\cite{Hakimi1962} provided constructive algorithms.  
Since verifying the Erd{\H o}s–Gallai inequalities and executing the Havel–Hakimi algorithm both run in polynomial time with the length of the degree sequences, the degree sequence problem for graphs is in  P.  Hypergraphs generalize graphs by allowing edges to contain arbitrary positive number of vertices~\cite{Berge1989}.  A hypergraph is called \emph{$t$-uniform} if each edge contains exactly $t$ vertices.  The hypergraphic degree sequence problem (HDS) asks, given a sequence of integers $D=(d_1,\dots,d_n)$ and an integer $t\ge 3$, whether there exists a simple (i.e., containing no multiple hyperedges) $t$-uniform hypergraph with degree sequence $D$.  A necessary condition is that $\sum_{i=1}^n d_i = t|E|$ for some integer $|E|$, by the generalized handshaking lemma, but unlike the simple graph case, no simple characterization is known.

Deza et al.~\cite{Deza2019} showed that  to decide if a sequence has a $3$-uniform  hypergraph realization is NP-complete.   In the dense regime, recent work by Li and Mikl\'os~\cite{LiMiklos2023} gives explicit degree thresholds $c_1 n^2$ and $c_2 n^2$ (with $0<c_1<c_2<\frac{1}{2}$) such that any sufficiently long degree sequence with all degrees between these bounds is always realizable by a $3$-uniform hypergraph - assuming that the sum is divisible by $3$. These results were extended to a dichotomy theorem by Logsdon \emph{et al.} \cite{lmmz2024}. The dichotomy theorem says that for any $c_2 \in (0,1)$ there exists a $c_1^*(c_2)$, such that for all $c_1^*(c_2)<c_1 \le c_2$, any degree sequence of length of sufficiently large  $n$ with all degrees between $c_1{n-1\choose 2}$ and $c_2{n-1\choose 2}$ has a 3-uniform hypergraph realization if and only if the sum of the degrees can be divided by $3$. On the other hand, for any $c_1 < c_1^*(c_2)$, it is NP-complete to decide if a degree sequence of length $n$ with degrees between $c_1{n-1\choose 2}$ and $c_2{n-1\choose 2}$ has a $3$-uniform hypergraph realization. Note that this is not a strict dichotomy theorem for the case $c_1 = c_1^*(c_2)$ is not resolved. Also note that for any $0<k<1$, $k{n-1\choose 2} = \Theta(n^2)$. 
 More generally, a complete $t$-uniform hypergraph on $n$ vertices has $\binom{n}{t}$ edges, so each vertex has degree $\binom{n-1}{t-1} = \Theta(n^{t-1})$.  Thus, a dense $t$-uniform hypergraph (for example, in the Erd{\H o}s-R\'enyi model with any fixed $p>0$) typically has minimum degree $\Omega(n^{t-1})$.

In this work, we study the opposite regime, that is, \emph{sparse} degree sequences, where the maximum degree is $o(n^{t-1})$.  Such sequences are natural in many applications: for example, hypergraphs arising from real-world data (e.g.\ social, biological, or informational networks) often have few edges per vertex.
Sparse hypergraphs form a central object of study in modern combinatorics and extremal graph theory.
A classical family of sparse hypergraphs arises from \emph{linear} designs, in which any two edges intersect in at most one vertex.
The finite projective plane of order $q$ is a canonical example: it is a $3$-uniform linear hypergraph with $n = q^2 + q + 1$ vertices and degree $q + 1 = \Theta(n^{1/2})$ \cite{Dembowski1968}.
Steiner systems $S(t,k,n)$ are highly regular $k$-uniform hypergraphs in which every $t$-subset of vertices is contained in exactly one block. It is well-known that the degrees in Steiner Triple Systems (in $S(2,3,n)$ systems) is $\frac{n-1}{2}$.
A natural probabilistic source of sparse examples is the Erd\H{o}s--R\'enyi random $t$-uniform hypergraph $H^{(t)}(n,p)$, where each $t$-set appears independently with probability $p$.
The expected degree is
\[
\mathbb{E}[d(v)] = p \binom{n-1}{t-1} = p \Theta( n^{t-1})=o(n^{t-1})
\]
for $p=p(n)=o(1)$. Particularly, when $p =o( n^{-(t-2)})$, we obtain $\mathbb{E}[d(v)] = o(n)$.
Sparse random hypergraphs played a central role in developing e.g., the powerful hypergraph container method, see for example~\cite{SaxtonThomason2015,BaloghMorrisSamotij2015}.
Hypergraphs of bounded or constant degree are trivially sparse yet may occur structurally rich.
High-dimensional expanders~\cite{Lubotzky2018} and simplicial complexes arising from Ramanujan-type constructions~\cite{Gromov2010} are examples of constant-degree hypergraphs that exhibit strong pseudorandomness and spectral expansion.
These constructions show that even highly sparse hypergraphs can maintain global combinatorial regularity and expansion properties.

One of the most challenging, widely unsolved problem is the Brown-Erd\H os-S\'os conjecture. Let $f_t(n,v,e)$ denote the maximum size of a $t$-uniform hypergraph on $n$ vertices such that no $v$ vertices span $e$ edges. They showed    
~\cite{BrownErdosSos1973,BrownErdosSos1973a} that $f_t(n,e(t-k)+k,e)=\Theta(n^k)$ but if we forbid $e$ edges on just one more vertex, this drops in magnitude, that is $f_t(n,e(t-k)+k+1,e)=o(n^k)$. Maybe the most famous solved special case ($t=3$, $e=3$, $k=2$) is the celebrated $(6,3)$-theorem of Ruzsa and Szemer\'edi \cite{RSZ} stating that $f_3(n,6,3)=o(n^2)$. 
For a beautiful, extremely important construction of sequences with no $k$-term arithmetic progression -- known to be sublinear by the famous theorem of Szemer\'edi -- see Behrend ~\cite{Behrend1946}.

These examples demonstrate that sparse hypergraphs can exhibit rich and highly constrained structure, underscoring the importance of understanding degree realizability in sparse regimes. They also motivate a systematic study of hypergraphicality under explicit sparsity constraints.
To this end, we consider a parametric version of the sparse $t$-uniform hypergraph degree sequence problem. For fixed $t \ge 3$ and constants $0 \le \alpha' \le \alpha < t- 1$, the problem $\text{sparse-}t\text{-uni-HDS}_{\alpha',\alpha}$ asks whether a given degree sequence admits a $t$-uniform hypergraph realization in which every degree lies between $n^{\alpha'}$ and $6n^{\alpha}$, where $n$ is the length of the sequence.

Our main result is a sharp dichotomy theorem establishing the exact boundary between polynomial-time solvable and NP-complete cases of this problem. Specifically, for every fixed $t \ge 3$, there exists a critical threshold
$$
\alpha^* := \frac{t(\alpha - 1) + 1}{t - 1},
$$
such that $\text{sparse-}t\text{-uni-HDS}_{\alpha',\alpha}$ is NP-complete whenever $\alpha' \le \alpha^*$, and solvable in linear time whenever $\alpha' > \alpha^*$.

This classification unifies and extends previous results obtained for dense hypergraphs, demonstrating that the computational complexity depends sharply on the interplay between the minimum and maximum degree exponents in the sparse degree sequence cases. 
It further reveals that even extremely sparse hypergraphs -- those with maximum degree $o(n)$ but at least $\Omega(n^{\frac{t-1}{t}})$ -- can retain the full NP-completeness of the general $t$-uniform realization problem. By contrast, when the degree range decreases beyond 
the threshold, realizability reduces to checking a simple divisibility condition, hence falling into linear time. This dichotomy links algorithmic complexity with structural sparsity in hypergraphs, providing a unified framework that bridges dense and sparse regimes and clarifies the computational landscape of hypergraphic degree sequences.

The remainder of the paper is organized as follows.  In Section~\ref{sec:prelim}, we give formal definitions and preliminaries.  In Section~\ref{sec:np-general}, we prove that the $t$-uniform HDS problem is NP-complete for each fixed $t\ge3$.  In Section~\ref{sec:np-sparse}, we prove the dichotomy theorem.  We conclude in Section~\ref{sec:conclusion} with a summary and discussion.

\section{Preliminaries}\label{sec:prelim}
A \emph{hypergraph} $H=(V,E)$ consists of a vertex set $V$ and an edge set $E\subseteq 2^V \setminus \{\emptyset\}$, where each edge $e\in E$ is a subset of $V$.  A hypergraph is \emph{$t$-uniform} if every edge contains exactly $t$ vertices.  For a vertex $v\in V$, the \emph{degree} $d(v)$ is the number of edges containing $v$.  The \emph{degree sequence} of $H$ is the nonnegative integer sequence of vertex degrees (usually written in nonincreasing order).  Given a sequence of integers $D=(d_1,d_2,\dots,d_n)$ and an integer $t$, the \emph{$t$-Uniform Hypergraphic Degree Sequence} (t-HDS) problem asks whether there exists a $t$-uniform hypergraph on the vertex set $\{v_1,v_2,\dots,v_n\}$ whose degree sequence is exactly $D$. We will work with vertex-labelled hypergraphs and indexed degree sequences, and throughout the paper, $d_i$ denotes the degree of vertex $v_i$.  We call such a sequence \emph{hypergraphic} if a realizing hypergraph exists. We will denote the degree of vertex $v$ in hypergraph $H$ by $d_H(v)$. When it is clear in which hypergraph $v$ is in, $H$ will be omitted from the subscript.

A degree sequence $D$ is \emph{regular} if all the degrees are the same in $D$. A degree sequence $D$ is \emph{almost regular} if there exists a $k$ such that each degree in $D$ is either $k$ or $k+1$. Observe that any regular degree sequence is also almost regular. We also call a subset of vertices regular (respectively, almost regular) if their degrees are regular (respectively, almost regular).

A basic necessary condition for hypergraphicality is that the sum of the degrees be divisible by $t$: since each edge contributes $1$ to the degree of each of its $t$ vertices, we have $\sum_{i=1}^n d_i = t\,|E|$ for a $t$-uniform hypergraph with $|E|$ edges.
Other than this divisibility condition, no simple criteria are known for $t$-uniform sequences.
 The t-HDS decision problem is clearly in NP: given a candidate hypergraph, one can verify the degrees in polynomial time.  We will show NP-hardness via reductions.

We use standard asymptotic notation.  A function $f(n)$ is $o(n)$ if $f(n)/n\to0$ as $n\to\infty$.  In particular, saying that a degree sequence has maximum degree $o(n)$ means each degree grows sublinearly in the number of vertices.  Finally, by a \emph{polynomial-time reduction} from a decision problem $A$ to $B$, we mean a polynomial-time computable transformation of instances of $A$ into instances of $B$ that preserves yes/no answers.

\section{NP-Completeness for General $t$-Uniform HDS}\label{sec:np-general}
We first show that for any fixed $t\ge3$, the $t$-uniform HDS problem is NP-complete.  
For $t = 3$, Deza et al.~\cite{Deza2019} proved that $3$-uniform HDS is NP-complete.  We extend this to all $t>3$ by a simple reduction.

\begin{theorem}\label{theo:general-t-uni}
For any fixed integer $t\ge3$, determining whether a given sequence of nonnegative integers is realizable as a $t$-uniform hypergraph degree sequence is NP-complete.
\end{theorem}

\begin{proof}
Hardness follows by induction from the known $3$-uniform case~\cite{Deza2019}.  Given an instance $D_t=(d_1,d_2,\dots,d_n)$ of $t$-uniform HDS, we construct in polynomial time a sequence $D_{t+1}$ for $(t+1)$-uniform HDS as follows.  Introduce a new vertex $v_{n+1}$ and define
\[
D_{t+1}=(d_1,\dots,d_n,d_{n+1}),
\quad\text{where}\quad d_{n+1} \;=\; \frac{d_1 + \cdots + d_n}{t}.
\]
Since $\sum_{i=1}^n d_i$ must be divisible by $t$ for $D_t$ to be $t$-hypergraphic, $d_{n+1}$ is an integer.  We claim $D_t$ has a $t$-uniform realization if and only if $D_{t+1}$ has a $(t+1)$-uniform realization.

Clearly, $t|\sum_{i=1}^{n}d_i$ iff $(t+1)|\sum_{i=1}^{n+1}d_i$.

($\Rightarrow$) Suppose $H_t$ is a $t$-uniform hypergraph realizing $D_t$.  Add the new vertex $v_{n+1}$ to every edge of $H_t$.  Since each of the original edges of $H_t$ contained $t$ vertices, adding the extra vertex yields a $(t+1)$-uniform hypergraph $H_{t+1}$ on $n+1$ vertices.  In $H_{t+1}$, the degree of vertex $v_{n+1}$ is exactly the number of edges of $H_t$, which is
$$
|E(H_t)| \;=\; \frac{1}{t}\sum_{i=1}^n d_i \;=\; d_{n+1},
$$
matching the construction of $D_{t+1}$.  Every original vertex $v_i$, $i\le n$ keeps its degree, since each $t$-uniform hyperedge incident to $v_i$ becomes a $t+1$-uniform hyperegde by adding $v_{n+1}$ to its vertex set. Hence $H_{t+1}$ realizes $D_{t+1}$.

($\Leftarrow$) Conversely, suppose there is a $(t+1)$-uniform hypergraph $H_{t+1}$ realizing $D_{t+1}$, $d(v_i)=d_i$ for $1\le i\le n+1$.  We argue that every edge of $H_{t+1}$ must contain the vertex $v_{n+1}$. Indeed, the total number of edges in $H_{t+1}$ is

$$\displaystyle{|E|= \frac{\sum_{i=1}^{n+1} d_i}{t+1}=\frac{\sum_{i=1}^{n} d_i}{t}=d_{n+1}.}$$

Removing $v_{n+1}$ from every edge yields a $t$-uniform hypergraph $H_t$ on vertices $v_1,\dots,v_n$ whose degree sequence is exactly $D_t$.  (It is straightforward to verify that the resulting hyperedges remain distinct and cover all prescribed degrees.)  Hence, $D_t$ is realizable if and only if $D_{t+1}$ is so.

The base case, that is, $3$-uniform HDS is NP-complete~\cite{Deza2019}.
The construction from $D_t$ to $D_{t+1}$ is polynomial-time computable.
By iterating this argument starting with the $3$-uniform HDS problem up to any fixed $t$, this reduction shows that $t$-uniform HDS is NP-hard by induction. Clearly the problem is in NP, so it is NP-complete.
Thus, we conclude that HDS is NP-complete for every fixed $t\ge3$.
\end{proof}

\section{Dichotomy theorem for the sparse $t$-uniform hypergraphicality}\label{sec:np-sparse}

First, we define the parameterized sparse $t$-uniform hypergraphicality problem.
\begin{definition}
    For any $t\ge 3$, and $0\le \alpha' \le \alpha <t-1$, the input of the $sparse-t-uni-HDS_{\alpha',\alpha}$ problem (sparse $t$-uniform hypergraph degree sequence problem with parameters $\alpha'$ and $\alpha$) has degree sequences of length $n$ as inputs, for which the degrees in the degree sequence are between $n^{\alpha'}$ and $6n^{\alpha}$ and asks the decision question if there exists a $t$-uniform hypergraph with the prescribed degrees.
\end{definition}
\begin{remark}
    The constant $6$ is needed to prove a sharp dichotomy theorem. Without the constant $6$, only a less strict dichotomy theorem could be proved that proves NP-completeness for $\alpha'$ to be strictly less than $ \frac{t(\alpha-1)+1}{t-1}$. The $6$ is not an absolute lower bound constant. It could be replaced with a function of $t$ that has upper bound $6$.
\end{remark}

We are going to prove a dichotomy theorem on this parametric decision problem saying that this decision problem is either NP-complete or can be solved in linear time.

First, we prove a lemma that we will use in both direction of the dichotomy theorem.
\begin{lemma}\label{lem:almost-regular-part}
Let $H = (V,E)$ be an arbitrary $t$-uniform hypergraph, and let $U\subseteq V$. Then there exists a $t$-uniform hypergraph  $H'$ with the following properties.
\begin{enumerate}[(i)]
\item For all $v\in V\setminus U$, $d_H(v) = d_{H'}(v)$,
\item $\sum_{v\in U} d_{H}(u) =\sum_{v\in U}d_{H'}(u)$,
\item the degrees in $U(H')$ are almost regular.
\end{enumerate}
\end{lemma}
\begin{proof}
We are going to construct a finite series of hypergraphs
$$
H = H_0, H_1, \ldots, H_k,
$$
where $H_k$ has the prescribed properties.
Starting with $l= 0$, let $H_l$ be the current constructed hypergraph in the sequence, and let $\bar{d}$ be the average degree of the vertices in $U(H_l)$. If the degrees of $U(H_l)$ is not almost regular, then there exists a $v_i\in U$ such that $d_{H_l}(v_i) > \lfloor\bar{d}\rfloor$ and there exists a $v_j \in U$ such that $d_{H_l}(v_j) <\lceil\bar{d}\rceil$. Find the largest such $v_i$ and the smallest such $v_j$. Then either $d_{H_l}(v_i) > \lceil\bar{d}\rceil$ or $d_{H_l}(v_j) < \lfloor\bar{d}\rfloor$ (or both) since $U$ is not almost regular in $H_l$. Also, there exists an $e\in E(H_l)$ such that $v_i \in e$, $v_j\notin e$ and $(e\setminus\{v_i\})\cup \{v_j\} \notin E(H_l)$. Construct $H_{l+1}$ from $H_l$ by removing $e$ and adding $(e\setminus\{v_i\})\cup \{v_j\} \notin E(H_l)$. $H_{l+1}$ is a simple $t$-uniform hypergraph in which for all $v\in V\setminus U$, $d_{H_l}(v) = d_{H}(v)$, and $\sum_{v\in U} d_{H}(u) =\sum_{v\in U}d_{H'}(u)$. Further,
$$
\sum_{v\in U} \min\{|d_{H_{l+1}}(v) - \lfloor\bar{d}\rfloor|,|d_{H_{l+1}}(v) - \lceil\bar{d}\rceil|\} < \sum_{v\in U} \min\{|d_{H_{l}}(v) - \lfloor\bar{d}\rfloor|,|d_{H_{l}}(v) - \lceil\bar{d}\rceil|\}.
$$
Since $f(G):=\sum_{v\in U} \min\{|d_{G}(v) - \lfloor\bar{d}\rfloor|,|d_{G}(v) - \lceil\bar{d}\rceil|\} =0$ if and only if $U(G)$ is almost-regular and $f(G)$ is an integer valued function, in finite number of steps we will obtain a hypergraph with the prescribed properties.
\end{proof}
{\it Remark:} The operation that removes edge $e$ and add  edge $(e\setminus \{v_i\})\cup\{v_j\}$ is called \emph{hinge-flip} operation, see for example, \cite{LiMiklos2023}.  

We need a technical lemma for the hardness (NP-completeness) direction.
\begin{lemma}\label{lem:almost-regular-gadget}
Given integers $t\ge 3$, $s>0$, $d>0$, and  $m>0$, such that ${m\choose t-1} \ge d$. Define a technical constant
$$
k := m\left\lceil\frac{sd(t-1)}{m}\right\rceil-sd(t-1).
$$
 Then the degree sequence consisting of $s$ degrees of $d$, $k$ degrees of $\left\lfloor\frac{sd(t-1)}{m}\right\rfloor$ and $m-k$ degrees of $\left\lceil\frac{sd(t-1)}{m}\right\rceil$ has a $t$-uniform hypergraph realization $H$, such that any edge in $H$ is incident with exactly one vertex with degree $d$. 
\end{lemma}
\begin{proof}
Let $V$ be the set of vertices with prescribed degree $d$ and let $U$ be the set of vertices with other prescribed degrees.
Construct a $\tilde{H}$ $t$-uniform hypergraph by adding arbitrarily $d$ different hyperedges to each vertex in $V$ such that these hyperedges incident with one vertex in $V$ and $t-1$ vertices in $U$. This is possible since $|U| =m$ and ${m\choose t-1}\ge d$. The average degree in $U(H)$ is $\frac{sd(t-1)}{m}$. According to Lemma~\ref{lem:almost-regular-part}, there exists a hypergraph with almost regular degrees in $U$. Since in $\tilde{H}$, every hyperedge is incident with $t-1$ vertices in $U$, and the hinge-flip operations in the proof of Lemma~\ref{lem:almost-regular-part} does not change this property, the hypergraph constructed by following the proof of Lemma~\ref{lem:almost-regular-part} will be a hypergraph with the prescribed degrees, and any edge in it is incident with exactly one vertex with degree $d$. 

\end{proof}

\noindent Now we are ready to prove the hardness direction of the dichotomy theorem.

\begin{theorem}[The hardness part of the dichotomy theorem]\label{theo:main}
For any fixed $t\ge3$, and parameters $(\alpha',\alpha)$ with $\alpha \in \left[\frac{t-1}{t},t-1\right)$ and $\alpha' \le \frac{t(\alpha-1)+1}{t-1}$, the $sparse-t-uni-HDS_{\alpha',\alpha}$ is NP-complete.
\end{theorem}

\begin{proof}
The proof is by a polynomial-time reduction from the general $t$-HDS problem, which is NP-complete, according to Theorem~\ref{theo:general-t-uni}.  Let $D=(d_1,\dots,d_m)$ be an arbitrary degree sequence of length $m$ (the original instance).  
We might assume that $m$ is bigger than a (later defined) threshold $m_0$. If $D$ has length shorther than $m_0$, then amalgamate $D$ by a sequence of $0$ degrees. Observe that the amalgamated degree sequence has a $t$-uniform hypergraph realization if and only if $D$ has one.
Let $n:= \left\lceil\left(\frac{6m}{t}\right)^{\frac{t-1}{\alpha}}\right\rceil$.
It follows that $m\in\left[\frac{t(n-1)^{\frac{\alpha}{t-1}}}{6},\frac{tn^{\frac{\alpha}{t-1}}}{6}\right]$.
We will construct a new degree sequence $D'$ of length $n$ whose maximum degree is at most $6n^{\alpha}$, its minimum degree is $\left\lceil n^{\frac{t(\alpha-1)+1}{t-1}}\right\rceil$, 
and such that $D'$ is realizable if and only if $D$ is realizable.
We obtain $D'$ by constructing a gadget hypergraph. The gadget hypergraph has such a rigid structure that forces any realization of $D'$ to decompose into a fixed part and an arbitrary realization of the original instance.

Partition the $n$ vertices into three sets:
\[
V_L,\quad V_N,\quad V_S,
\]
where $|V_N|=m$, $|V_L| = m$, and $|V_S| = n - 2m$ is the remainder (so $|V_S| = \Theta(n)$ since 
$m\le \frac{tn^{\frac{\alpha}{t-1}}}{6} = o(n)$
for $\alpha<t-1$).  We construct a $t$-uniform gadget hypergraph $G$ on these vertices with the following edges:
\emph{i)} all possible $t$-subsets of $V_L$ (so $V_L$ forms a complete $t$-uniform hypergraph).
\emph{ii)} all $t$-subsets in $V_L\cup V_n$ that each of them includes at least one vertex from both vertex sets,
\emph{iii)} and for each vertex $s\in V_S$,  $d= \left\lceil n^{\frac{t(\alpha-1)+1}{t-1}}\right\rceil$ hyperedges that contains $s$ together with $t-1$ vertices from $V_L$; these are chosen so that the extra degree contributions to $V_L$ from the $V_S$ vertices are spread almost regularly. Since $m\ge\frac{t(n-1)^{\frac{\alpha}{t-1}}}{6}$, ${m\choose t-1}\ge d$ for sufficiently large $m$, and thus such a configuration of hyperedges exists according to Lemma~\ref{lem:almost-regular-gadget}, with $s=n-2m$.
Here we define $m_0$ such that for all $m\ge m_0$, ${m\choose t-1}\ge d$.
Note that $G$ is not unique, however, we need only its degree sequence $D_G$, which is unique. 

Compute the degree sequence $D_G$ of this gadget $G$.  Every vertex in $V_L$ has degree coming from edges in $V_L$, edges between $V_L$ and $V_N$, and  $\left\lfloor\frac{(n-2m)d(t-1)}{m}\right\rfloor$ or  $\left\lceil\frac{(n-2m)d(t-1)}{m}\right\rceil$ edges incident with one vertex in $V_S$.  One checks easily that the degree of each $v\in V_L$ is
either
$$
{2m-1\choose t-1} + \left\lceil\frac{(n-2m)d(t-1)}{m}\right\rceil
$$
or
$$
{2m-1\choose t-1} + \left\lfloor\frac{(n-2m)d(t-1)}{m}\right\rfloor.
$$
Since $m\in\left[\frac{t(n-1)^{\frac{\alpha}{t-1}}}{6},\frac{tn^{\frac{\alpha}{t-1}}}{6}\right]$ and $d = \left\lceil n^{\frac{t(\alpha-1)+1}{t-1}}\right\rceil$, this yields degree
$$
\deg_G(v) \le {2\left\lfloor\frac{t n^{\frac{\alpha}{t-1}}}{6}\right\rfloor-1\choose t-1}+(n-2m)(t-1)\frac{\left\lceil n^{\frac{t(\alpha-1)+1}{t-1}}\right\rceil}{\frac{t(n-1)^{\frac{\alpha}{t-1}}} {6}} <
$$
$$
<\frac{\left(\frac{t}{3}\right)^{t-1}n^\alpha}{(t-1)!}+\frac{6(t-1)}{t}\frac{n-2m}{(n-2m)^{\frac{\alpha}{t-1}}}\left\lceil n^{\frac{t(\alpha-1)+1}{t-1}}\right\rceil< \frac{\left(\frac{t}{3}\right)^{t-1}n^\alpha}{(t-1)!}+\frac{6(t-1)}{t}n^{\alpha}<6n^{\alpha}.
$$
 Vertices in $V_N$ get edges only from the mixed $V_L$--$V_N$ family; each vertex $v\in V_N$ has degree ${2m-1\choose t-1}-{m-1\choose t-1} < 6n^{\alpha}$.  Each vertex in $V_S$ has exactly degree $d =\left\lceil n^{\frac{t(\alpha-1)+1}{t-1}}\right\rceil$ by construction. It is easy to see that all degrees are larger than or equal to $n^{\alpha'}$. Therefore, the degrees are in the prescribed interval.

Now let $D'$ be the sequence obtained by taking $D_G$ and adding the entries of the original sequence $D$ to the degrees of the vertices in $V_N$ (i.e.\ for $v_i \in V_N$ we set $d'_i = \deg_H(v_i)+d_i$ and all other degrees do not change).
The result is a degree sequence of length $n$.
Since for each $i \le m$, $d_i \le {m-1\choose t-1}$, the maximum degree in $D'$ is still less than $6n^{\alpha}$.
We claim that $D'$ has a realization if and only if $D$ does.

If $D$ is realizable, let $H_N$ be a $t$-uniform hypergraph on $V_N$ realizing $D$.  Combine $H_N$ with the gadget hypergraph $G$: that is, take the union of their edge sets (treating them as hypergraphs on the disjoint vertex sets $V_N$ and $V_L\cup V_S$, respectively). Since the induced sub-hypergraph on $V_N$ in $G$ is empty, the union will be a simple hypergraph. In this union, each vertex in $V_N$ now has its degree 
as the degree in the gadget hypergraph $G$ plus the degree $d_i$ 
(the degree from $H_N$), yielding exactly $d'_i$.  Each vertex in $V_L\cup V_S$ keeps the degree it had in $G$, which was designed to satisfy $D_G$.  Thus, this union hypergraph realizes $D'$.

Conversely, suppose $D'$ is realizable by some $t$-uniform hypergraph $H'$.  By analyzing the degree constraints, we show  that any realization of $D'$ must decompose into some realization $G'$ of $D_G$ plus a $t$-uniform hypergraph realizing $D$. Intuitively, the gadget is so rigid (complete on $V_L$ and fully connected to $V_N$ from $V_L$) that any realization of $D'$ must be the disjoint union of a realization $G'$ of $D_G$ and $D$.

Let $D_L$ be the maximal degrees in $D'$, let $D_S$ be the minimal degrees in $D'$ and let $D_N$ be the set of remaining degrees. We call the corresponding vertex sets by $V_L$, $V_S$ and $V_N$.

The sum of the degrees in $V_L$ is
$$
|V_L|{|V_L| +|V_N|-1 \choose t-1}+ |V_S|d(t-1) = m{2m-1\choose t-1} +(n-2m)d(t-1).
$$
Let $x$ denote the number of hyperedges in $H'$ that are incident to only vertices from $V_L$. Let $n_i$ denote the number of hyperedges in $H'$ that are incident with $i$ vertices in $V_N$ and with $t-i$ vertices in $V_L$. Let $s_j$ denote the number of hyperedges in $H'$ that are incident with $j$ vertices from $V_S$ and with $t-j$ vertices in $V_L$. Finally, let $q_{i,j}$ denote the hyperedges in $H'$ that are incident with $i$ vertices from $V_N$, with $j$ vertices in $V_S$ and $t-i-j$ vertices in $V_L$. Then we get that
\begin{eqnarray}
&&m{2m-1\choose t-1} +(n-2m)d(t-1) = \nonumber\\
&&= tx+\sum_{i=1}^{t-1}n_i(t-i) + \sum_{j=1}^{t-1}s_j(t-j)+\sum_{j=1}^{t-1}\sum_{i=1}^{t-j-1}q_{i,j}(t-i-j).  \label{eq:degree-sum-constraint}
\end{eqnarray}
We know that $x\le {m\choose t}$ and $n_i \le {m\choose i}{m\choose t-i}$, since these are the maximum number of such kind of hyperedges. Furthermore, 
$$
\sum_{j=1}^{t-1} s_j j + \sum_{j=1}^{t-1}\sum_{i=1}^{t-j-1} q_{i,j}j \le (n-2m)d,
$$
for $(n-2m)d$ is the sum of the degrees in $V_S$. With this constraint,
$$
\sum_{j=1}^{t-1}s_j(t-j)+\sum_{j=1}^{t-1}\sum_{i=1}^{t-j-1}q_{i,j}(t-i-j)
$$
is maximal when $s_1 = (n-2m)d$, and for all $j>1$ $s_j=0$ and for all $i,j$, $q_{i,j}=0$. 
Indeed, the largest sum of degrees in $V_L$ can be obtained by using hyperedges that have as many vertices in $V_L$ as possible, as few vertices in $V_S$ as possible and no vertices in $V_N$.
In that case,
$$
\sum_{j=1}^{t-1}s_j(t-j)+\sum_{j=1}^{t-1}\sum_{i=1}^{t-j-1}q_{i,j}(t-i-j) = (n-2m)d(t-1).
$$
Therefore, equation~\ref{eq:degree-sum-constraint} is satisfied if and only if $x={m\choose t}$, for all $i$, $n_i = {m\choose i}{m\choose t-i}$, $s_1 = (n-2m)d$, for all $j>1$, $s_j = 0$ and for all $i,j$, $q_{i,j} = 0$.
That is, after removing edges incident to any vertex in $V_L$, we get a realization of $D$ in $V_N$.

Since the reduction from $D$ to $D'$ is computable in polynomial time, this shows that $t$-uniform HDS remains NP-hard even when restricted to sequences with $\max_i d_i \le 6n^{\alpha}$ and $\min_i d_i \ge n^{\frac{t(\alpha-1)+1}{t-1}}$.  Being in NP is unchanged.  This completes the proof.
\end{proof}

Now we are going to prove the tractability direction of the dichotomy theorem. It follows directly from the Theorem~\ref{theo:fully-graphic-regime}. We need a lemma to prove Theorem~\ref{theo:fully-graphic-regime}. That lemma has been proved in previous papers in slightly different form, here we give a proof for completeness. The lemma is about the \emph{least balanced degree sequences} that we define below.
\begin{definition}
Let $n$ be a positive integer, and let $0\le \delta \le \Delta$. Further, let  $\Sigma$ be a non-negative integer with $\delta n\le \Sigma\le \Delta n$. The \emph{least balanced degree sequence}, $LBDS(n,\delta,\Delta,\Sigma)$ is the degree sequence of length $n$ that contains only degrees $\delta$, $\Delta$, at most one degree $d$ with $\delta < d < \Delta$, and whose sum of degrees is $\Sigma$. 
\end{definition}
It is easy to see that the $LBDS(n,\delta,\Delta,\Sigma)$ exists, and when it is written in non-increasing order, it is unique. We prove the following lemma.
\begin{lemma}\label{lem:lbds}
    Let $D$ be a degree sequence of length $n$, with minimum degree $\delta$ maximum degree $\Delta$ and degree sum $\Sigma$. Then it holds that $\delta n\le \Sigma\le \Delta n$. If the $LBDS(n,\delta,\Delta,\Sigma)$ has a $t$-uniform hypergraph realization, then $D$ also has a $t$-uniform hypergraph realization.
\end{lemma}
\begin{proof}
    First, we give a finite series of perturbations that transforms $D$ to the $LBDS(n,\delta,\Delta,\Sigma)$. Then we use the inverse of this series to construct a realization of $D$ from a realization of the $LBDS(n,\delta,\Delta,\Sigma)$.
    The finite series of perturbations goes in the following way. If $D$ is not the $LBDS(n,\delta,\Delta,\Sigma)$, there are at least two degrees which are larger than $\delta$ and smaller than $\Delta$. Find a smallest and a largest such degree, denote it by $d_i$ and $d_j$, $d_i\le d_j$, and create a new degree sequence by decreasing $d_i$ by $1$ and increasing $d_j$ by $1$. By iterating these steps, either $d_i$ becomes $\delta$ or $d_j$ becomes $\Delta$ (or both). Then the number of vertices which are neither $\delta$ nor $\Delta$ decreases at least by $1$. Therefore, this procedure reaches the $LBDS(n,\delta,\Delta,\Sigma)$ in finite number of steps. Let the so-obtained series of degree sequences be
    $$
    D = D_0, D_1, \ldots, D_k = LBDS(n,\delta,\Delta,\Sigma).
    $$
    Let $H_k$ be a $t$-uniform realization of the $LBDS(n,\delta,\Delta,\Sigma)$. Starting with $l=k$, while $l\neq 0$, do the following. Let $d_i$ and $d_j$ be the degrees for which decreasing $d_i$ by $1$ and increasing $d_j$ by $1$ $D_{l-1}$ was transformed to $D_l$. Let the corresponding vertices be denoted by $v_i$ and $v_j$. Take the $t$-uniform realization of $D_l$, and find a hyperedge $e$ such that $v_j\in e$, $v_i\notin e$ and $(e\setminus \{v_j\})\cup\{v_i\}\notin E(H_l)$. Such a hyperedge exists for $d_i < d_j$. We just obtained a realization of $D_{l-1}$. Iterate this procedure till we get a realization of $D_0 = D$.
\end{proof}

\noindent We are ready to prove the tractability direction of the dichotomy theorem.

\begin{theorem}\label{theo:fully-graphic-regime}
    For any $t\ge 3$, and any $\alpha'\le\alpha<t-1$ with $\alpha' > \frac{t(\alpha-1)+1}{t-1}$, there exists a computable $n_0$ such that for all $n\ge n_0$, any degree sequence of length $n$ with minimum degree at least $n^{\alpha'}$ and maximum degree $6n^\alpha$ has a $t$-uniform hypergraph realization if and only in the sum of the degrees is divisible by $t$. 
\end{theorem}
\begin{proof}
By the Lemma~\ref{lem:lbds}, it suffices to prove that for all sufficiently large $n$, any $\alpha' > \frac{t(\alpha-1)+1}{t-1}$ and for all $\Sigma$,  such that $t|\Sigma$ and $n\left\lceil n^{\alpha'}\right\rceil\le \Sigma \le n\left\lfloor 6n^{\alpha}\right\rfloor$, the $LBDS(n, \left\lceil n^{\alpha'}\right\rceil,\left\lfloor 6n^{\alpha}\right\rfloor, \Sigma)$ has a $t$-uniform hypergraph realization. 
Let $D$ be the $LBDS(n, \left\lceil n^{\alpha'}\right\rceil,\left\lfloor 6n^{\alpha}\right\rfloor, \Sigma)$, and let $k$ denote the number of degrees of $6n^{\alpha}$ in $D$. There might or might not be a degree of $d$ with $n^{\alpha'} < d < 6n^{\alpha}$. We are going to construct a $t$-uniform hypergraph realization of $D$, denoted by $H = (V,E)$. Let $V_L$ be the set of vertices with prescribed degree $6n^{\alpha}$, and let $V_S$ be the set of vertices with prescribed degree $n^{\alpha'}$. We consider the following cases, in each case we consider the presence or absence of degree $d$.
\begin{enumerate}[(i)]
\item $k\le \frac{n^{\alpha'+1-\alpha}}{12(t-1)}$. Since $\alpha'\le \alpha$, $k\le \frac{n}{12(t-1)}$. Therefore, ${n-k-2\choose t-1} > 6n^{\alpha}$ for sufficiently large $n$ since $\alpha < t-1$. Further, the corresponding threshold $n_0$ is easy to compute. Therefore, for each vertex $v$ with degree $6n^{\alpha}$ or $d$ we can select $6n^{\alpha}$ or $d$ $t-1$-tuples from $V_S$. Together with $v$, they will be $t$-uniform hyperedges. Let $H'$ be the $t$-uniform hypergraph obtained by adding these hyperedges by selecting the above-mentioned $t-1$ tuples. In $H'$ every vertex $v\in V_L$ and the possible vertex with prescribed degree $d$ have their prescribed degree. 
According to Lemma~\ref{lem:almost-regular-part}, there exists a hypergraph $H"$ in which $V_S$ is almost regular, and every vertex in $V_L$ and the possible vertex with prescribed degree $d$ have the prescribed degrees.
Since $(t-1)(k6n^{\alpha}+d) < (n-k-1)n^{\alpha'}$ for sufficiently large $n$, the average degree in $V_S$ in $H"$ is smaller than $n^{\alpha'}$. Further, the corresponding threshold $n_0$ is easy to compute. Since ${|V_S|-1\choose t-1} > n^{\alpha'}$, there exists an almost regular $t$-uniform hypergaph $H'''$ on vertices in $V_S$ such that the union of the hyperedges in $H"$ and $H'''$ is a realization of the $D$. Note that no hyperedge in $H"$ is incident with vertices only in $V_S$, while $H'''$ all hyperedges incident with vertices only in $V_S$, thus, the union of $H"$ and $H'''$ is a simple $t$-uniform hypergraph. 

\item $\frac{n^{\alpha'+1-\alpha}}{12(t-1)} < k \le \frac{t-2}{t-1}n-1$. Since $\alpha' > \frac{t(\alpha-1)+1}{t-1}$, $\alpha'+1-\alpha > \frac{\alpha}{t-1}$. Therefore, for sufficiently large $n$, ${k-1\choose t-1} > 6n^{\alpha}$. Further, the corresponding threshold $n_0$ is easy to compute. Therefore, for sufficiently large $n$, there exists a hypergraph on $V_L$, such that each vertex has a degree $6n^{\alpha}$ except at most $t-1$ that has degree $6n^{\alpha}-1$ (so the sum of the degrees is divisible by $t$).
Indeed, it is an almost regular $t$-uniform hypergraph whose existence follows from Lemma~\ref{lem:almost-regular-part}.
Since for sufficiently large $n$ (with an easy-to-compute threshold $n_0$), $(t-1+d)(t-1) < |V_S|$, it is possible to select $t-1+d$ instances of disjoint $t-1$-tuples from $V_S$ as incident vertices of $t$-uniform hyperedges that adjust the degrees in $V_L$ to the prescribed degree and also satisfy the prescribed degree $d$ of a vertex, if such vertex exists. Further, ${|V_S|-1\choose t-1} > n^{\alpha'}$, therefore, the so-obtained hypergraph can be extended to a realization of $D$. Indeed, this extension is an almost uniform hypergraph on $V_S$, whose existence follows from Lemma~\ref{lem:almost-regular-part}.
\item $k> \frac{t-2}{t-1}n+1$. In this case, $k6n^{\alpha}>(t-1)((n-k-1)n^{\alpha'}+d)$. Therefore, we can select $n^{\alpha'}$ or $d$ $t-1$-tuples from $V_L$ for each vertex with prescribed to degree $n^{\alpha'}$ or $d$ to satisfy their degrees, and the average degree in $V_L$ will be less than $6n^{\alpha}$. 
Applying Lemma~\ref{lem:almost-regular-part} with $U = V_L$, there exists a hypergraph $H'$ for which the degrees in $V_L$ is almost-regular. Since ${|V_L|-1\choose t-1} > 6n^{\alpha}$, we can extend $H'$ to a realization of $D$ since this extension is again an almost regular hypergraph on $V_L$.
\end{enumerate}
We considered all possible cases regarding the value of $k$, thus the proof is complete.
\end{proof}

Now we can prove the tractability direction of the dichotomy theorem.
\begin{theorem}
    For any $t\ge 3$, and for all $0\le\alpha'\le \alpha<t-1$ with $\alpha'>\frac{t(\alpha-1)+1}{t-1}$, the $sparse-t-uni-HDS_{\alpha',\alpha}$ problem can be solved in linear time.
\end{theorem}
\begin{proof}
    From Theorem~\ref{theo:fully-graphic-regime}, there exists a computable $n_0$, such that for all $n> n_0$, any input degree sequence of length $n$ has a $t$-uniform hypergraph realization if and only if the sum of the degrees is divisible by $t$. Computing $n_0$ does not depend on the input, so it can be done in constant time. 
    When $n<n_0$, the input size is constant, and even a brute-force algorithm can check in constant time whether the input degree sequence is hypergraphic. 
    When $n\ge n_0$, deciding that the sum of the degrees can be divided by (a constant value) $t$ can be done in linear time.
\end{proof}

\section{Conclusion}\label{sec:conclusion}

We have established a full dichotomy theorem for the parameterized sparse $t$-uniform hypergraph degree sequence problem, $\mathrm{sparse}\text{-}t\text{-}\mathrm{uni}\text{-}\mathrm{HDS}_{\alpha',\alpha}$. Our result precisely identifies the threshold $\alpha' = \frac{t(\alpha - 1) + 1}{t - 1}$ that separates NP-complete and linear time solvable instances. This classification provides a complete picture of the computational complexity landscape of hypergraph degree sequence realizability, bridging the gap between previously known dense and sparse regimes.

Beyond its immediate implications for hypergraph theory, the dichotomy reveals a deeper structural phenomenon: the transition from algorithmic intractability to tractability is governed not by absolute degree magnitudes but by the asymptotic ratio between minimum and maximum degree exponents. This insight unifies and extends multiple strands of prior research on graphic and hypergraphic degree sequences, offering a general conceptual framework for reasoning about sparsity and complexity in combinatorial realization problems.

Future work may explore analogous phase transitions in more general combinatorial settings, including relations among minimum, maximum degrees and sum of the degrees. For example, the results presented in this paper might be compared with the work of Aldesari and Greenhill \cite{ag2019}. They provided an asymptotic enumeration formula for the number of $t$-uniform hypergraph realizations of degree sequences. Their formula works for degree sequences with maximum degree $d$ and sum of degrees $\Sigma$ when
\begin{equation}
t^4 d^3 = o(\Sigma).\label{eq:ag-condition}
\end{equation}
It is easy to see that a necessary condition to satisfy Equation~\ref{eq:ag-condition} is that $d$ should be $o(n^{\frac{1}{2}})$. 
A natural open question if the results of Aldesari and Greenhill can be extended to degree sequences with maximum degree $o(n^{\frac{t-1}{t}})$. Also, it is an interesting question if other dichotomy theorems could be obtained with other parametrization, for example, parameterizing the degree sequences with the maximum degree and sum of the degrees.

\section*{Acknowledgement} The work was supported by the Hungarian NKFIH grant K132696 and by the European Union project RRF2.3.1-21-2022-00006 within the framework of Health Safety National Laboratory Grant no RRF-2.3.1-21-2022-00006.

\end{document}